\documentclass[11pt]{article}
\usepackage[english]{babel}
\usepackage{amssymb,amsmath,amsthm}
\usepackage[latin2]{inputenc}
\usepackage{amsfonts, verbatim, amscd}
\usepackage{pstricks}

\reversemarginpar
\def\ol{\overline}
\def\wt{\widetilde}
\def\wh{\widehat}

\newtheorem{theorem}{Theorem}
\newtheorem{proposition}[theorem]{Proposition}
\newtheorem{lemma}[theorem]{Lemma}
\newtheorem{remark}[theorem]{Remark}
\def\acknowledgement{\par\addvspace{17pt}\small\rmfamily
\trivlist\if!\ackname!\item[]\else
\item[\hskip\labelsep
{\bfseries\ackname}]\fi}

\def\eps{\varepsilon}
\def\phi{\varphi}
\def\bs{\backslash}

\def\C{\mathbb{C}}
\def\R{\mathbb{R}}

\def\cF{\mathcal{F}}

\def\cO{{\cal O}}

\newcommand{\parc}[2]{\frac{\partial {#1}}{\partial {#2}}}

\begin{document}

\centerline{\bf CLOSED HOLOMORPHIC 1-FORMS WITHOUT ZEROS}
\centerline{\bf ON STEIN MANIFOLDS}
\vskip 5mm
\centerline{by}
\vskip 5mm
\centerline{IRENA MAJCEN}
\bigskip

\it
\centerline{IMFM, University of Ljubljana}
\centerline{Ljubljana, Slovenia} 
\rm
\vskip 15mm

\section{Introduction}\label{intro}
A \emph{Stein} manifold is a complex manifold biholomorphic to a closed complex submanifold of a complex Euclidean space $\C ^N$. Stein manifolds of complex dimension one are precisely open Riemann surfaces.

Suppose that $X$ is a Stein manifold. In every class of the cohomology group $H^1(X,\C)$ there is a closed holomorphic 1-form (Theorem 1 in \cite{Ser}, see also \cite[p.\ 160]{Car} and Theorem 2 in \cite[p.\ 208]{Ste}). The existence of a holomorphic function $f\colon X\to\C$ without critical points, proved by Forstneri\v{c} in \cite{For1}, implies that in the zero class there is a closed holomorphic 1-form without zeros, namely $df$. (For the case when $X$ is an open Riemann surface see also \cite{GunNar}.) Our goal in this paper is to show that a closed holomorphic 1-form without zeros can be chosen in every cohomology class. (See also \cite{KusSai} for open Riemann surfaces.)

\begin{theorem}\label{glavni}
Let $X$ be a Stein manifold. Every cohomology class in the cohomology group $H^1(X,\C)$ is represented by a closed holomorphic 1-form without zeros.
\end{theorem}

Note that Theorem \ref{glavni} is not true for an arbitrary complex manifold $X$. For example, if $X$ is a compact Riemann surface of genus $g$ then by the Riemann-Roch theorem each closed holomorphic 1-form has precisely $2g-2$ zeros.

Denote by $\cO (X)$ the algebra of all holomorphic functions on $X$. A compact set $K\subset X$ is said to be \emph{$\cO (X)$-convex} if for any point $x\in X\bs K$ there exists $f\in \cO(X)$ such that $|f(x)|>\max_K |f|$. 

Theorem \ref{glavni} is a consequence of the following result.

\begin{theorem}\label{applic}
Let $K\subset X$ be a compact, $\cO(X)$-convex subset and $\theta$ a closed 1-form on $X$. Let $\omega$ be a closed holomorphic 1-form on a neighborhood $U$ of $K$, with $\omega| _x\ne 0$ for all $x\in K$, such that $\int_C \omega =\int _C \theta$ for each closed curve $C\subset U$. Then there exists  a closed holomorphic 1-form  $\wt\omega$ without zeros on $X$ such that $[\wt\omega]=[\theta]\in H^1(X,\C)$, and there exists a holomorphic injective mapping $h$ in a neighborhood of $K$, close to the identity, such that $\wt\omega = h^* \omega$ near $K$.
\end{theorem}

A closed nowhere vanishing holomorphic 1-form $\omega$ on $X$ defines a holomorphic foliation $\cF$ of codimension 1 with the tangent bundle $T\cF=\ker \omega$. Theorem \ref{applic} then gives an approximation theorem for holomorphic nonsingular hypersurface foliations, defined by nonvanishing closed holomorphic 1-forms. The approximation is done by a global foliation $\widetilde \cF$ on $X$ such that $\cF$ and $\widetilde \cF$ are conjugate to each other on a neighborhood of $K$.
\bigskip

Necessary and sufficient conditions for a compact real manifold without boun\-dary to admit a representative without zeros in each class of the de Rham cohomology group have been studied in \cite{Lat}.

\section{Proof of the main theorems}
We only need to prove Theorem \ref{applic}. It can be assumed that $\theta$ in Theorem \ref{applic} is a holomorphic 1-form on $X$ (Theorem 1 in \cite{Ser}). We show that there exists a closed holomorphic 1-form $\omega$ such that $\omega|_x\ne 0$ for all $x\in X$ and
$$\int_{C_i} \theta = \int _{C_i}\omega$$
for every closed curve $C_i\subset X$ where $\{C_1,C_2,\ldots\}$ is a basis of $H_1(X,\R)$.

\medskip
Choose a number $\wh c\in \R$. Fix a smooth strongly plurisubharmonic Morse exhaustion function $\rho \colon X\to \R$ such that $\rho < \wh c$ in $K$ but $\rho > \wh c$ in $U^c$ (\cite{Hor}, Theorem 5.1.6.).
The sublevel sets $\{\rho\le c\}$ are compact and $\cO (X)$-convex for all $c\in \R$, and for every regular value $c$ of $\rho$, $\{\rho < c\}$ is a smooth strongly pseudoconvex domain. We can replace $K$ by $\{\rho\le \wh c\}$ and in the following subsections the notations $K$ and $U$ will be used for other purposes.

Suppose that a closed holomorphic nonvanishing 1-form $\omega$ is given in a neighborhood of $\{\rho \le c\}$, where $c$ is a regular value of $\rho$, and suppose that $\omega$ has the same periods as $\theta$ on each closed curve $C_i\subset \{\rho\le c\}$. The main step in the proof is a construction of a closed holomorphic nowhere vanishing 1-form $\wt \omega$ on a neighborhood of $\{ \rho \le \wt c\}$, which has the same periods as $\theta$, where $\wt c>c$ is also a regular value of $\rho$. We will obtain a global solution using a limit process.

Our construction requires two distinct arguments depending on whether the interval $[c,\wt c]$ contains a critical value of $\rho$ or not and is similar to the scheme in Section 6.~in \cite{For1}. 

We choose a distance function $d$ on $X$ induced by a smooth Riemannian metric on $TX$. Suppose that $\alpha\colon V\to X$ is a holomorphic mapping on an open set $V$ in $X$. We denote $||\alpha-id||_V=\sup_{x\in V} d(\alpha(x),x)$.

We first cite a result which will be needed several times in the proof.

\begin{theorem}[Theorem 4.1 in \cite{For1}]\label{split}
Let $A$ and $B$ be compact sets in a complex manifold $X$ such that $A\cup B$ has a basis of Stein neighborhoods in $X$ and $\ol{A\bs B} \cap \ol{B\bs A}=\emptyset$. Given an open set $\wt C\subset X$ containing $C:=A\cap B$ there exist open sets $A'\supset A$, $B'\supset B$, $C'\supset C$ with $C'\subset A'\cap B'\subset \wt C$, satisfying the following. For every $\eta>0$ there is $\eps_\eta >0$ such that for each injective holomorphic map $\gamma \colon \wt C\to X$ with $||\gamma - id||_{\wt C} < \eps_\eta$ there exist injective holomorphic maps $\alpha \colon A'\to X$, $\beta\colon B'\to X$ satisfying
$$\gamma =\beta \circ\alpha^{-1} \textrm{ on } C',\quad ||\alpha - id||_{A'} < \eta,\quad ||\beta - id||_{B'} < \eta.$$
\end{theorem}

\subsection{The noncritical case}\label{noncritical}


Let $A\subset \wt A \subset X$ be compact sets in $X$. If there exists a smooth strongly plurisubharmonic function $\rho$ in an open set $\Omega \supset \ol{\wt A \bs A}$ which has no critical points on $\Omega$ and satisfies
$$	A\cap \Omega     = \{x\in \Omega \colon \rho(x)\le 0\},\quad   
        \wt A\cap \Omega = \{x\in \Omega \colon \rho(x)\le 1\},$$
we call $\wt A$ a {\it noncritical strongly pseudoconvex extension} of $A$. The set $A_t=A\cup \{\rho \le t\} \subset X$ is a smooth strongly pseudoconvex domain in $X$ for each $t\in [0,1]$ and the family smoothly increases from $A=A_0$ to $\wt A=A_1$. A manifold $X$ is said to be a {\it noncritical strongly pseudoconvex extension} of $A$ if 
there is a smooth exhaustion function $\rho\colon X \to \R$
such that $A=\{\rho\le 0\}$ and $\rho$ is strongly plurisubharmonic and without critical points on $\{\rho\ge 0\}=X\bs {\rm int A}$.

With this notations the noncritical case can be rephrased as follows.

\begin{proposition}\label{prop2}
Let $\theta$ be a closed 1-form on a Stein manifold $X$. Suppose that $\widetilde A\subset X$ is a noncritical strongly pseudoconvex extension of $A\subset \widetilde A$. Let $\omega$ be a closed holomorphic 1-form on a neighborhood of $A$ with $\omega|_x\ne 0$ for all $x\in A$ such that $\int _C\omega = \int_C \theta$ for all closed curves $C\subset A$. Choose $\eps>0$. There exists a closed holomorphic 1-form $\widetilde\omega$ on a neighborhood of $\widetilde A$, with $\widetilde \omega|_x \ne 0$ for all $x\in \widetilde A$, such that $\int _C\wt\omega = \int_C \theta$ for all closed curves $C\subset \wt A$, and there exists a holomorphic injective mapping $\alpha$ in a neighborhood $V$ of $A$ such that $\widetilde\omega=\alpha^*\omega$ near $A$ and $||\alpha-id||_V < \eps$.
\end{proposition}

\begin{proof}
We first introduce the notion of a convex bump. We shall use the same kind of bumps as in \cite{For1} and \cite{For2}. Denote the coordinates on $\C^n$ by $z=(z_1, \ldots, z_n)$ where $z_j=x_j+iy_j$ and let
$$P=\{z\in \C^n\colon |x_j|< 1,\ |y_j|<1,\ j=1,2, \ldots, n\}$$
denote the open unit cube. Set $P'=\{z\in P\colon y_n=0\}.$ We say that a compact set $B\subset X$ is a \emph{convex bump} on a compact set $A\subset X$ if there exist an open set $U\subset X$ containing $B$, a biholomorphic map $\phi \colon U\to P$ onto $P\subset \C^n$ and smooth strongly concave functions $h, \wt h \colon P'\to [-a,a]$ for some $a<1$ such that $h\le \wt h$ and $h=\wt h$ near the boundary of $P'$ and
$$\phi(A\cap U)=\{z\in P\colon y_n\le h(z_1, \ldots, z_{n-1}, x_n)\},$$
$$\phi((A\cup B)\cap U)=\{z\in P\colon y_n\le \wt h(z_1, \ldots, z_{n-1}, x_n)\}.$$

Assume that $A\subset \wt A$ is a strongly pseudoconvex extension. Applying Lemma 12.3.~in \cite{HenLei1} we find finitely many compact strongly pseudoconvex domains $A=A_0 \subset A_1 \subset \ldots \subset A_{k_0} = \wt A$ in $X$ such that $A_{k+1}=A_k\cup B_k$ for every $k=0,1,\ldots, k_0-1$ where $B_k$ denotes a convex bump on $A_k$ as defined above. (For similar constructions see also \cite{For1}, \cite{ForPre1}, \cite{ForPre2}, \cite{ForPre3}, \cite{Gro}, \cite{HenLei2}.)

The construction breaks into $k_0$ steps of the same kind. In the $k$-th step we show how to obtain a holomorphic 1-form $\omega_{k+1}$ with prescribed periods on a neighborhood of $A_{k+1}$ if we are given $\omega_k$ on a neighborhood of $A_k$ which satisfies
$$\omega_k=\alpha_{k-1}^* \omega_{k-1} \textrm{ in } A'_{k-1} \quad \textrm{and}\quad ||\alpha_{k-1} - Id||_{A'_{k-1}}\le \frac{\eps}{2^{k}},$$
where $A'_{k-1}$ is a neighborhood of $A_{k-1}$.

\begin{figure}[htbp]
\centering
\begin{pspicture}(-4,-2.3)(4,2)
\pscustom[linestyle=none, fillstyle=solid, fillcolor=gray]{
	\psarc(0,-6){6}{63.256}{116.744}
	\psline(-2.7,-0.642)(-2.7,-2)(2.7,-2)(2.7,-0.642)
	}
\pscustom[linestyle=none, fillstyle=vlines]{
	\pscurve(-1.553,-0.204)(-2,-0.343)(-2.4,-0.7)(-2,-1)(2,-1)(2.4,-0.7)(2,-0.343)(1.553,-0.204)(0,0)(-1.553,-0.204)
}
\psarcn(0,-6){6}{120}{60}
\psframe[linestyle=dashed](-3,-2.3)(3,1.9)
\uput[r](3,1){$U$}
\psframe[dimen=inner](-2.7,-2)(2.7,1.6)
\uput[d](-2,1.6){$\wt L$}
\pscurve(-1.553,-0.204)(-0.7,0.2)(0,0.4)(0.7,0.2)(1.553,-0.204)
\uput*[u](2,-2){$\wt K$}
\psecurve(-1.553,-0.204)(-2,-0.343)(-2.4,-0.7)(-2,-1)(2,-1)(2.4,-0.7)(2,-0.343)(1.553,-0.204)
\rput*(0,-0.5){$C_k$}
\end{pspicture}
\caption{Convex bump}
\end{figure}
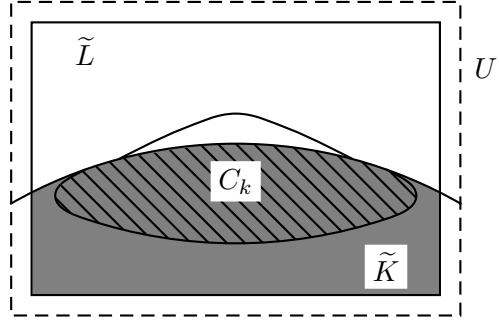

Choose a neighborhood $U$ of $B_k$ as in the definition of a convex bump. With the notations as in the definition choose $c\in (a,1)$ sufficiently close to 1 such that the compact support of $\widetilde h -h$ is contained in $cP'$. Let $L:=c\overline{P}\subset \C^n$ and $\widetilde L :=\phi^{-1}(L)\subset U$. By letting $c$ tend to 1, we may assume that $B_k\subset \widetilde L$. Set $\widetilde K=A_k\cap \widetilde L$ and $K=\phi(\widetilde K)\subset P$. Note that $K$ is convex and hence $\wt K$ has a simply connected neighborhood $V$.

There is a noncritical holomorphic function $f$ in a neighborhood $V$ of $\wt K$ such that $\omega_k=df$ on $V$. Choose $\eps'>0$. By Proposition 3.3.~in \cite{For1} we obtain a holomorphic submersion $g\colon V'\to \C$ in an open set $V'\supset \wt L$ which approximates $f$ uniformly in a neighborhood of $\wt K$.

Denote $C_k=A_k\cap B_k$. By Lemma 5.1.~in \cite{For1} there exist constants $\eps_0>0$, $M>0$ and an open set $W\subset X$ with $C_k\subset W \subset V$ satisfying the following property. Given $\eps'\in (0,\eps_0)$ and a holomorphic submersion $g\colon W\to \C$ with $\sup_{x\in W} |f(x)-g(x)|<\eps'$ there is an injective holomorphic map $\gamma \colon W\to X$ satisfying $f=g\circ \gamma$ on $W$ and $||\gamma - id||_{W} < M\cdot \eps '$. By the above argument such a holomorphic submersion $g$ exists, hence we have the decomposition $f=g\circ \gamma$ on $W$. But $\gamma$ can also be decomposed.

Set $\eta=\frac{\eps}{2^{k+1}}$ and $\eps'=\min\{ \frac{\eps_{\eta}}{M}, \frac{\eps_0}{2}\}$. By the above construction the mapping $\gamma \colon W \to X$ can be chosen to satisfy $||\gamma-id||_{W}< \eps _\eta$ and can be then split using Theorem \ref{split} as
$$\gamma = \beta_k \circ \alpha_k^{-1}\quad \textrm{with}\quad ||\alpha_k - id||_{A_k'} < \frac{\eps}{2^{k+1}},\quad ||\beta_k - id||_{B_k'} < \frac{\eps}{2^{k+1}}.$$

Define a holomorphic 1-form $\omega_{k+1}$ by
$$\omega _{k+1} = \left\{
\begin{array}{cc}
\alpha_k^* \omega_k & \textrm{on } A_k',\\
d(g\circ \beta_k) & \textrm{on } B_k'.
\end{array}\right.$$
Clearly $\omega_{k+1}$ is well defined and without zeros. The mapping $\alpha_k$ provides a homotopy of curves $C$ and $\alpha_k(C)$ for any closed curve $C\subset A_0$, thus the relation $\omega_{k+1} = \alpha_k^* \omega_k$ implies that 1-forms $\omega_{k+1}$ and $\omega_k$ have the same periods (by Stokes' theorem).

Repeating this step $k_0$ times we obtain $\wt\omega=\omega_{k_0}$ with the required properties in a neighborhood of $\wt A=A_{k_0}$ satisfying $\wt\omega=h^* \omega$ on $A_0$ where 
$h=\alpha_1 \circ \alpha_2 \circ \ldots \circ \alpha_{k_0}$. By construction we have
\begin{eqnarray*}
||h-id||_{A_0} &= &||\alpha _1 \circ \ldots \circ \alpha_{k_0} -\alpha_2\circ\ldots\circ\alpha_{k_0}+ \ldots + \alpha_{k_0} -id||_{A_0}\\
&\le & \sum_{k=1}^{k_0-1} ||\alpha _k \circ \ldots \circ \alpha_{k_0} -\alpha_{k+1}\circ\ldots\circ\alpha_{k_0}||_{A_0} + ||\alpha_{k_0} -id||_{A_0}\\
&< & \frac{\eps}{2} + \ldots + \frac{\eps}{2^{k_0-1}}+ \frac{\eps}{2^{k_0}} < \eps,
\end{eqnarray*}
since we may assume that $\alpha_2\circ\ldots\circ\alpha_{k_0}(A_0)\subset A_0'$, where $A_0'\supset A_0$ is a domain of $\alpha_1$.
\qed\end{proof}

With this we have finished the proof of the noncritical case. In the following subsection we treat the critical case.

\subsection{The critical case}\label{critical}
Let $p$ be a critical point of $\rho$ and let $k$ denote the Morse index of $p$. If $k=0$, $\rho$ has a local minimum at $p$; as $c$ passes $\rho(p)$, a new connected component appears in $\{\rho<c\}$ (see Lemma 2.3.~in \cite{HenLei2}). Hence $\omega$ can be trivially extended by taking a differential of any noncritical holomorphic function near $p$.

From now on we assume $k\ge 1$. There is no loss of generality in assuming $\rho(p)=0$. Choose $c_0\in (0,1)$ such that $p$ is the only critical point of $\rho$ in $[-c_0,3c_0]$. In what follows we explain how to obtain a closed nowhere vanishing holomorphic 1-form $\wt \omega$ with prescribed periods in a neighborhood of $\{\rho \le +c_0\}$ provided that a required 1-form $\omega$ has already been constructed on a neighborhood of $\{\rho \le -c_0\}$.
\bigskip

Denote by $P\subset \C^n$ the open unit polydisc. Using Lemma 2.5 in \cite{HenLei2} we may assume that all the critical points of $\rho$ are \emph{nice}, meaning that there is a neighborhood $U\subset X$ of $p$ and a biholomorphic coordinate map $\phi\colon U\to P$, with $\phi(p)=0$, such that 
the function $\widetilde \rho(z) := \rho(\phi^{-1}(z))$ is given by
$$\widetilde \rho(z)
= \sum_{j=1}^n \mu_j y_j^2 - \sum_{j=1}^k x_j^2 + \sum _{j=k+1}^n x_j^2,$$
where $\mu_j\ge 1$ for all $j$ and $\mu _j >1$ when $1\le j\le k$. Denote the smallest of the numbers $\mu_1, \ldots, \mu_k$ by $\mu$.


Write $z=(z',z'')=(x'+iy', x''+iy'')$, where $z'\in\C^k$ and $z''\in \C^{n-k}$. The set $E_0\subset U$ defined by 
$$
	\phi(E_0) = \{(x'+iy',z'')\colon  y'=0,\ z''=0,\ |x'|^2\le c_0 \}
$$
is a $k$-dimensional core of a handle attached from the outside to 
$\{\rho\le -c_0\}$. Such cores were used in \cite{For1}. The handles introduced in \cite{HenLei2} could have been used as well, but the first aproach has been chosen for it will be used in the proof of a more general result in Section 3.

Let $c=(1-\frac{1}{\mu})^2 c_0$. By the noncritical case
we may assume that $\omega$ is given on $\{\rho < -c/2 \}$. Define a smaller $k$-dimensional handle $E\subset E_0$ by the condition
$$\phi(E) = \{(x'+iy',z'')\colon  y'=0,\ z''=0,\ |x'|^2\le c \}.$$
Note that $bE$ is a $(k-1)$-sphere contained in $\{\rho = -c\}$.

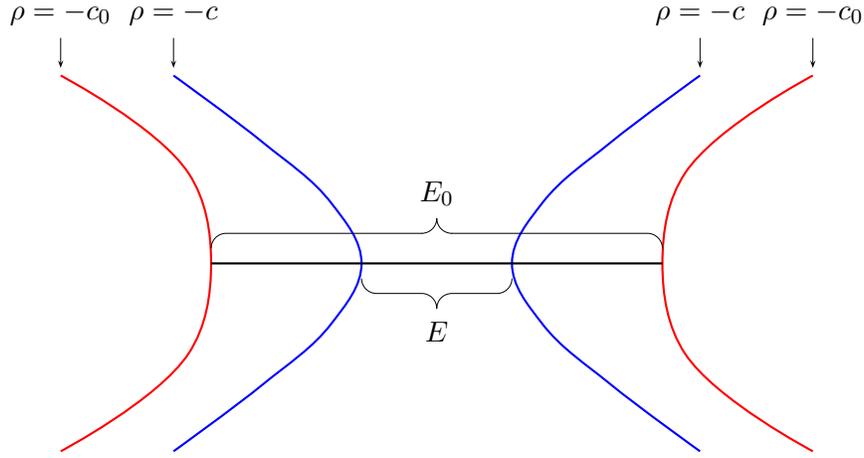
\begin{figure}[htb]
\centering
\begin{pspicture}(-5.5,-2.5)(5.5,3.5)
\pscurve[linecolor=red](-5,-2.5)(-3.3,-1.2)(-3,0)(-3.3,1.2)(-5,2.5)
\pscurve[linecolor=red](5,-2.5)(3.3,-1.2)(3,0)(3.3,1.2)(5,2.5)
\pscurve[linecolor=blue](-3.5,-2.5)(-2.3,-1.6)(-1.4,-0.8)(-1,0)(-1.4,0.8)(-2.3,1.6)(-3.5,2.5)
\pscurve[linecolor=blue](3.5,-2.5)(2.3,-1.6)(1.4,-0.8)(1,0)(1.4,0.8)(2.3,1.6)(3.5,2.5)
\psline(-3,0)(3,0)

\pscustom[linewidth=0.3pt]{
	\psarcn(-2.8,0.2){0.2}{180}{90}
	\psline(-2.8,0.4)(-0.2,0.4)
	\psarc(-0.2,0.6){0.2}{270}{360}
	\psarc(0.2,0.6){0.2}{180}{270}
	\psline(0.2,0.4)(2.8,0.4)
	\psarcn(2.8,0.2){0.2}{90}{0}
}

\uput[u](0,0.6){$E_0$}

\pscustom[linewidth=0.3pt]{
	\psarc(-0.8,-0.2){0.2}{180}{270}
	\psline(-0.8,-0.4)(-0.2,-0.4)
	\psarcn(-0.2,-0.6){0.2}{90}{0}
	\psarcn(0.2,-0.6){0.2}{180}{90}
	\psline(0.2,-0.4)(0.8,-0.4)
	\psarc(0.8,-0.2){0.2}{-90}{0}
}
\uput[d](0,-0.6){$E$}

\pscustom[linewidth=0.3pt]{
	\psline{->}(-5,3)(-5,2.6)\psline[liftpen=2]{->}(5,3)(5,2.6)
	\psline[liftpen=2]{->}(-3.5,3)(-3.5,2.6)\psline[liftpen=2]{->}(3.5,3)(3.5,2.6)
}
\uput[u](-5,3){$\rho=-c_0$}\uput[u](5,3){$\rho=-c_0$}
\uput[u](-3.5,3){$\rho=-c$}\uput[u](3.5,3){$\rho=-c$}

\end{pspicture}
\caption{Handles $E$ and $E_0$}
\end{figure}


When $k\ge 3$ we have $\pi_1(bE)=0$, hence $\int_\Gamma \theta=0$ for any closed curve $\Gamma\subset bE$. When $k=2$, $bE$ is a circle bounding the 2-disc $E$ and hence $\int_{bE} \omega = \int_{bE}\theta = \int _E d\theta =0$. In both cases we conclude that there exist a neighborhood $V$ of $bE$ in $\{\rho\le -\frac{c}{2}\}$ and a holomorphic submersion $f\colon V\to \C$ such that $\omega = df$ on $V$.


By Lemma 6.4.~in \cite{For1} there is a constant $c'\in (\frac12c,c)$ such that $f$ and its partial derivatives $\parc{f}{z_l}$ extend smoothly to $(\{\rho \le -c'\} \cap V)\cup E$ (without changing their values on $\{\rho \le -c'\}$) and the Jacobian matrix $J(\wt f)=(\parc{f_j}{z_l})$ of the extension $\wt f$ has complex rank $q$ at each point of $E$. Inspection of the proof shows that $c'$ can be chosen arbitrarily close to $c$, hence we may assume that
$$\phi^{-1} \Big(\{ (x'+iy',z'')\colon y'=0,\ z''=0,\ c'\le |x'|^2\le c\}\Big) \subset V.$$

Let $d(z,E)=\inf_{w\in E} d(z,w)$. The function $d^2(z,E)$ is strictly plurisubharmonic in a neighborhood of $E$ and provides a family of pseudoconvex neighborhoods of $E$.  Fix a compact pseudoconvex set $L$ such that $E\subset L \subset U$ and $L\cap \{\rho\le - c'\}\subset V$. Choose $c''\in (c',c)$ and denote $K=(L\cap\{\rho \le -c''\}) \cup E$. Lemma 6.6.~in \cite{For1} gives for every $\delta>0$ an open neighborhood $V'\subset X$ of the set $K$ and a holomorphic submersion $g\colon V'\to \C$ such that $|\wt f- g|_K<\delta$, $|d\wt f - dg|_E < \delta$. Here $|f|_K$ is the uniform norm of $f$ on $K$ and $|df|_E$ is the norm of its differential on $E$, measured in a fixed Hermitean metric on $TX$.


As in the noncritical case we shall use $g$ to define a holomorphic 1-form on a larger domain. The main difference is that this larger domain is not a sublevel set with some convex bumps but a sublevel set of a new strongly plurisubharmonic function $\tau$.

\begin{lemma}[Lemma 6.7.~in \cite{For1}]
There exists a smooth strongly plurisubharmonic function $\tau$ on $\{\rho < 3c_0\} \subset X$ which has no critical values in $(0, 3c_0)\subset \R$ and satisfies
\begin{enumerate}
\item[(i)] $\{\rho \le -c_0\} \cup E \subset \{\tau \le 0\} \subset \{\rho\le -c\} \cup E$ and
\item[(ii)] $\{\rho \le c_0\} \subset \{\tau \le 2c_0\} \subset \{\rho < 3c_0\}.$
\end{enumerate}
\end{lemma}



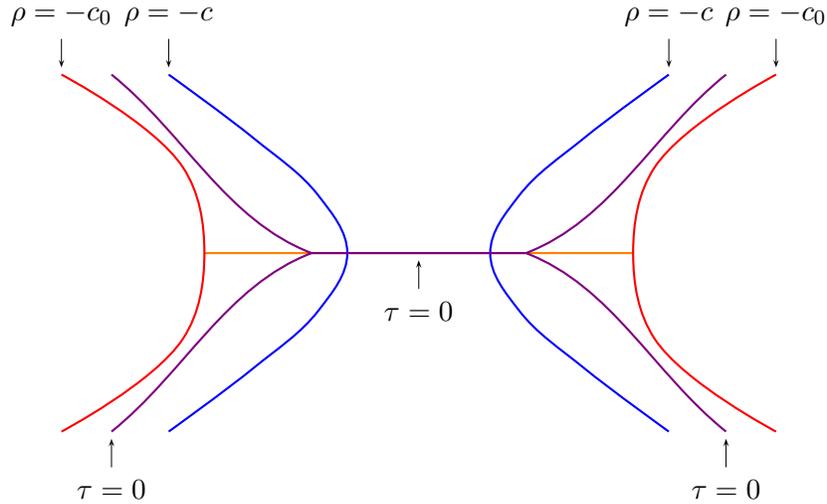
\begin{figure}[htb]
\centering
\psset{unit=0.95cm}
\begin{pspicture}(-5.5,-3.2)(5.5,3.3)
\psline[linecolor=orange](-3,0)(3,0)
\pscurve[linecolor=red](-5,-2.5)(-3.3,-1.2)(-3,0)(-3.3,1.2)(-5,2.5)
\pscurve[linecolor=red](5,-2.5)(3.3,-1.2)(3,0)(3.3,1.2)(5,2.5)
\pscurve[linecolor=blue](-3.5,-2.5)(-2.3,-1.6)(-1.4,-0.8)(-1,0)(-1.4,0.8)(-2.3,1.6)(-3.5,2.5)
\pscurve[linecolor=blue](3.5,-2.5)(2.3,-1.6)(1.4,-0.8)(1,0)(1.4,0.8)(2.3,1.6)(3.5,2.5)

\pscustom[linewidth=0.3pt]{
	\psline{->}(-5,3)(-5,2.6)\psline[liftpen=2]{->}(5,3)(5,2.6)
	\psline[liftpen=2]{->}(-3.5,3)(-3.5,2.6)\psline[liftpen=2]{->}(3.5,3)(3.5,2.6)
	\psline[liftpen=2]{->}(-4.3,-3)(-4.3,-2.6)\psline[liftpen=2]{->}(4.3,-3)(4.3,-2.6)
	\psline[liftpen=2]{->}(0,-0.5)(0,-0.1)
}
\uput[u](-5,3){$\rho=-c_0$}\uput[u](5,3){$\rho=-c_0$}
\uput[u](-3.5,3){$\rho=-c$}\uput[u](3.5,3){$\rho=-c$}

\psline[linecolor=violet](-1.5,0)(1.5,0)
\psecurve[linecolor=violet](-0.5,0)(-1.5,0)(-4.3,2.5)(-5,3)
\psecurve[linecolor=violet](-0.5,0)(-1.5,0)(-4.3,-2.5)(-5,-3)
\psecurve[linecolor=violet](0.5,0)(1.5,0)(4.3,2.5)(5,3)
\psecurve[linecolor=violet](0.5,0)(1.5,0)(4.3,-2.5)(5,-3)
\uput[d](-4.3,-3){$\tau=0$}
\uput[d](4.3,-3){$\tau=0$}
\uput[d](0,-0.5){$\tau=0$}
\end{pspicture}
\caption{The level set $\{\tau = 0\}$}
\end{figure}

Choose $c'''\in (c'',c)$. The set $\Omega=\{\rho < c'''\}\cup V'$ is a neighborhood of $\{\rho\le c\}\cup E$. Consider a family of sublevel sets $\{\tau \le t\}$ as $t$ increases from 0 to $2c_0$. The property (i) of $\tau$ implies that for sufficiently small $t>0$ we have $\{\tau \le t\} \subset \Omega$. Fix a number $t$ with such property.



\begin{figure}[htb]
\centering
\begin{pspicture}(-5.5,-3.2)(5.5,3.3)
\pscustom[fillstyle=solid, fillcolor=yellow, linestyle=none]{
	\psline(5,0)(5,2.5)(3.9,2.5)
	\psecurve(5.5,3.3)(3.9,2.5)(1.5,0.3)(0.82,0.185)(0,0.18)
	\psline(0.82,0.185)(0.8,0)(5,0)
}

\pscustom[fillstyle=solid, fillcolor=yellow, linestyle=none]{
	\psline(5,0)(5,-2.5)(3.9,-2.5)
	\psecurve(5.5,-3.3)(3.9,-2.5)(1.5,-0.3)(0.82,-0.185)(0,-0.18)
	\psline(0.82,-0.185)(0.8,0)(5,0)
}
\pscustom[fillstyle=solid, fillcolor=yellow, linestyle=none]{
	\psline(-5,0)(-5,2.5)(-3.9,2.5)
	\psecurve(-5.5,3.3)(-3.9,2.5)(-1.5,0.3)(-0.82,0.185)(0,0.18)
	\psline(-0.82,0.185)(-0.8,0)(-5,0)
}
\pscustom[fillstyle=solid, fillcolor=yellow, linestyle=none]{
	\psline(-5,0)(-5,-2.5)(-3.9,-2.5)
	\psecurve(-5.5,-3.3)(-3.9,-2.5)(-1.5,-0.3)(-0.82,-0.185)(0,-0.18)
	\psline(-0.82,-0.185)(-0.8,0)(-5,0)
}

\psline[linecolor=orange](-3,0)(3,0)
\pscurve[linecolor=red](-5,-2.5)(-3.3,-1.2)(-3,0)(-3.3,1.2)(-5,2.5)
\pscurve[linecolor=red](5,-2.5)(3.3,-1.2)(3,0)(3.3,1.2)(5,2.5)
\pscurve[linecolor=blue](-3.3,-2.5)(-2.1,-1.6)(-1.2,-0.8)(-0.8,0)(-1.2,0.8)(-2.1,1.6)(-3.3,2.5)
\pscurve[linecolor=blue](3.3,-2.5)(2.1,-1.6)(1.2,-0.8)(0.8,0)(1.2,0.8)(2.1,1.6)(3.3,2.5)

\pscustom[linewidth=0.3pt]{
	\psline{->}(-5,3)(-5,2.6)\psline[liftpen=2]{->}(5,3)(5,2.6)
	\psline[liftpen=2]{->}(-3.3,3)(-3.3,2.6)\psline[liftpen=2]{->}(3.3,-3)(3.3,-2.6)
}
\uput[u](-5,3){$\rho=-c_0$}\uput[u](5,3){$\rho=-c_0$}
\uput[u](-3,3){$\rho=-c''-\mu$}\uput[d](3.6,-3){$\rho=-c''-\mu$}

\psecurve[linecolor=violet](5.5,3.3)(3.9,2.5)(1.5,0.3)(0.7,0.19)(0,0.18)(-0.7,0.19)(-1.5,0.3)(-3.9,2.5)(-5.5,3.3)
\psecurve[linecolor=violet](5.5,-3.3)(3.9,-2.5)(1.5,-0.3)(0.7,-0.19)(0,-0.18)(-0.7,-0.19)(-1.5,-0.3)(-3.9,-2.5)(-5.5,-3.3)
\uput[u](3.9,3){$\tau=t$}\psline[linewidth=0.3pt]{->}(3.9,3)(3.9,2.6)
\uput[d](-3.9,-3){$\tau=t$}\psline[linewidth=0.3pt]{->}(-3.9,-3)(-3.9,-2.6)
\psframe(-2.6,-1.5)(2.6,1.5)
\uput[u](0,1.5){$L'$}

\pscustom[fillstyle=hlines, linestyle=none]{
	\psecurve(3.9,2.5)(2.6,1.2)(1.5,0.3)(0.7,0.19)(0,0.18)(-0.7,0.19)(-1.5,0.3)(-2.6,1.2)(-3.7,2.5)
	\psline[liftpen=1](-2.6,1.2)(-2.6,0)(2.6,0)(2.6,1.2)
}
\pscustom[fillstyle=hlines, linestyle=none]{
	\psecurve(-3.9,-2.5)(-2.6,-1.2)(-1.5,-0.3)(-0.7,-0.19)(0,-0.18)(0.7,-0.19)(1.5,-0.3)(2.6,-1.2)(3.9,-2.5)
	\psline[liftpen=1](2.6,-1.2)(2.6,0)(-2.6,0)(-2.6,-1.2)
}

\rput(0,2.5){$A$}
\pscustom[linewidth=0.3pt]{
	\psline(0.2,2.5)(2,2.5)\psline[liftpen=2]{->}(2,2.5)(3.5,1)
	\psline[liftpen=2](-0.2,2.5)(-2,2.5)\psline[liftpen=2]{->}(-2,2.5)(-3.5,1)
	\psline[liftpen=2]{->}(0,-2)(0,-0.1)
}
\uput[d](0,-2){$B$}

\end{pspicture}
\caption{Compact sets $A$ and $B$}
\end{figure}
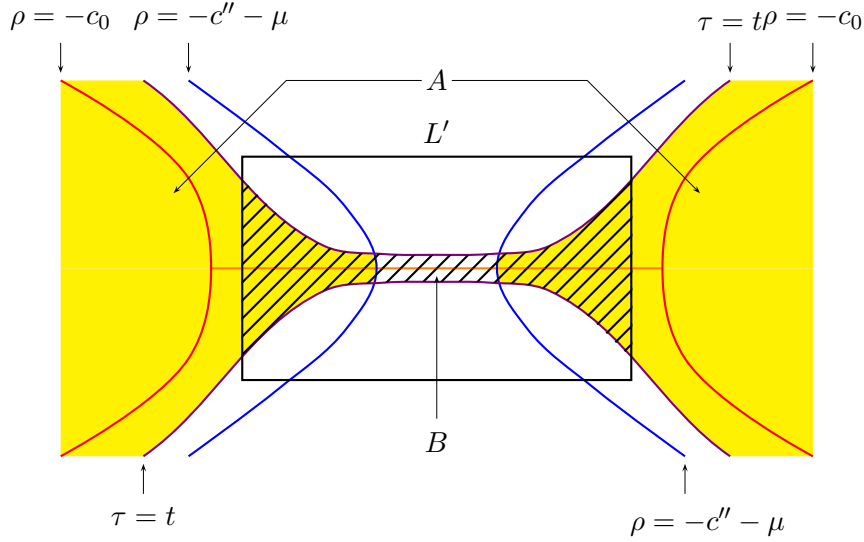


Choose a small number $\mu >0$. Let $L'$ be a compact pseudoconvex set such that $E\subset L'\subset \textrm{int} L$ and let $K'=(L'\cap\{\rho\le -c''-\mu\}) \cup E$. Finally, define
$$A=\{\tau \le t\} \cap \{\rho\le -c''-\mu\},\quad  B=L' \cap \{\tau \le t\}.$$
Note that $A\cup B=\{\tau \le t\}$ and $A\cap B= \{\rho\le -c''-\mu\} \cap \{\tau \le t\} \cap L'$, hence $W=\textrm{int } K$ is a neighborhood of $A\cap B$. Taking $\mu$ small enough, we have $\ol{A\bs B} \cap \ol{B\bs A} = \emptyset$. By Lemma 5.1. in \cite{For1} we obtain positive constants $\eps_0, M$ and an open neighborhood $W'\subset W$ of $A\cap B$ such that the following holds. For any $\eps' \in (0, \eps_0)$ and a holomorphic submersion $g\colon W\to \C$ with $\sup_{x\in W} |f(x)-g(x)|<\eps'$ there is an injective holomorphic map $\gamma \colon W' \to X$ satisfying
$$
f=g\circ\gamma \textrm{ on } W' \quad \textrm{and} \quad ||\gamma- id||_{W'}< M \cdot \eps'.
$$
The submersion $g$ is provided by the argument above for $\delta < \eps_0$. Additionally, to split $\gamma$ using Theorem \ref{split} and obtain an estimate $||\alpha-id||_{A'}<\frac{\eps}{2}$, we take $\eta=\frac{\eps}{2}$ and choose $\delta>0$ so small that $\delta \cdot 2M < \eps_\eta$. Note that $\wt C$ coresponds to $W'$. Thus Theorem \ref{split} gives injective holomorphic maps $\alpha$ and $\beta$ mapping from neighborhoods $A'$ of $A$, respectively $B'$ of $B$, into $X$ such that
$$\gamma=\beta\circ\alpha^{-1} \textrm{ on } C',\quad ||\alpha - id||_{A'} < \frac\eps2, \quad ||\beta - id||_{B'} < \frac\eps2.$$

Finally define
$$ \omega ' = \left\{
\begin{array}{cc}
\alpha^* \omega & \textrm{on } A'\\
d(g\circ \beta) & \textrm{on } B'.
\end{array}\right.$$
Note that $\omega'$ is given on a neighborhood of $\{\tau \le t\}$ and by the noncritical case (Proposition \ref{prop2}) we further obtain $\wt \omega$ on a neighborhood of $\{ \tau \le 2c_0\} \supset \{\rho \le c_0\}$ as desired.
\bigskip


It remains to construct $\wt \omega$ when $k=1$. There are two possibilities:
\begin{enumerate}
\item[(1)] the core of a handle $E_0$ connects two components of $\{\rho\le -c_0\}$,
\item[(2)] the core of a handle $E_0$ is a part of a new closed curve in the group $H_1(\{\rho \le c_0\}, \R)$.
\end{enumerate}
When the situation is as in (1), $\wt \omega$ can be constructed as above. The possibility (2) requires some additional effort. We can choose a new curve $C_0$ in such a way that $E_0\subset C_0$ and $C_0\bs E_0 \subset \{\rho\le -c_0\}$. The construction of $\wt\omega$ is similar as in the case $k\ge 2$, so to obtain $\wt\omega=h^*\omega$, the additional requirement is that $\int _{C_0} \wt \omega = \int _{C_0} \theta$. 

Since $bE$ consists of two points $a$ and $b$, we can choose disjoint neighborhoods $V_a$ and $V_b$ of $a$ and $b$ such that $\omega|_{V_a} = df_a$ and $\omega|_{V_b}=df_b$, where $f_a\colon V_a\to \C$ and $f_b\colon V_b\to \C$ are noncritical functions satisfying $f_b(b)-f_a(a)=\int_{C_0}\theta - \int_{C_0\bs E} \omega$. Denote $V=V_a\cup V_b$.

Applying the construction when $k\ge 2$ we get $\wt \omega = d(g\circ \beta)$ on a neighborhood of $E$ for some noncritical function $g$. The integral condition implies that $f\circ\alpha(b)-f\circ\alpha(a)+\int_{C_0\bs E}\alpha^*\omega=\int_{C_0}\theta$. In order to find such $\alpha$ and choose the appropriate $f$, we define a family of maps.

For each $t\in \Delta \subset \C$ let $f_t$ be a noncritical function such that
$$f_t \colon V\to \C,\quad f_t|_{V_a} = f_a \quad \textrm{ and } \quad f_t|_{V_b} =f_b+t.$$
We proceed as in the case $k\ge 2$. Using Lemma 6.4.~in \cite{For1} for a family of functions $f_t$ that depend holomorphically on parameter $t$, we obtain functions $\wt f_t$ on $(\{\rho \le c'\}\cap V)\cup E$ which smoothly extend $f_t$ and its partial derivatives $\parc{f_t}{z_l}$, and the Jacobian matrices $J(\wt f_t)$ have complex rank 1 at each point of $E$. We define $L$ and $K$ as before and use Lemma 6.6.~in \cite{For1} to obtain a family $g_t\colon V' \to \C$ such that
$$|\wt f_t-g_t|_K< \delta \quad \textrm{and}\quad |d\wt f_t-dg|_E<\delta.$$

Let $A$ and $B$ be defined as above. Lemma 5.1.~in \cite{For1} gives maps $\gamma_t$, holomorphic in $t$, such that $f_t=g_t\circ\gamma_t$ on a neighborhood of $K$. Define $\gamma(x,t):=(\gamma_t(x),t)$. By Theorem 4.1.~in \cite{For1} for sets $A\times\Delta$ and $B\times\Delta$ we obtain a splitting $\gamma=\beta\circ\alpha^{-1}$ where the mappings $\alpha$ and $\beta$ are of the same form as $\gamma$, that is $\alpha(x,t)=(\alpha_t(x),t)$ and $\beta(x,t)=(\beta_t(x),t)$. Thus we get a family of 1-forms defined by
$$\omega_t=\left\{
\begin{array}{cc}
\alpha_t^{*} \omega& \textrm{on } A',\\
d(g_t\circ \beta_t)& \textrm{on } B'.
\end{array}\right.$$

We now show that there is a number $t'\in \Delta$ such that $\int_{C_0} \omega _{t'} = \int _{C_0} \theta$. Define
$$\phi(t):=t, \quad \psi(t):=f_t(\alpha_t(b))-f_t(\alpha_t(a))+\int_{C_0\bs E}\alpha^*\omega - \int_{C_0}\theta.$$
Obviously $\psi$ and $\phi$ are holomorphic and for a fixed $s<1$, $s$ close to 1, we have $|\phi(t)-\psi(t)|$ small for $t\in b(s\Delta)$ and $|\phi(t)|=|t|=s$ for $t\in b(s\Delta)$, hence
$$|\phi(t)-\psi(t)|<|\phi(t)|$$
on $b(s\Delta)$. By Rouche's theorem $\phi$ and $\psi$ have the same number of zeros on $s\Delta$. Since $\phi(0)=0$, there is a $t'\in s\Delta$ such that $\phi(t')=0$. The form $\omega_{t'}$ satisfies all the requirements.

\subsection{Conclusion}

We now finish the proof of Theorem \ref{applic} by explaining the global scheme. Denote by $p_1, p_2, \ldots$ the critical points of $\rho$ in $\{\rho>\wh c\}$. Since each critical level set contains a unique critical point we may assume that
$$\rho(p_1)<\rho(p_2)<\rho(p_3)<\ldots\ .$$

We inductively choose a sequence of real numbers $\wt c_0<\wt c_1<\wt c_2<\ldots$ such that
$$\wt c_0 <\rho(p_1) < \wt c_1 < \rho (p_2) < \ldots \ .$$
Let $d_j= \rho(p_j)-\wt c_{j-1}$. By choosing $\wt c_{j-1}$ sufficiently close to $\rho(p_j)$ it can be achieved for each $j=1,2,\ldots$ that $\rho(p_j)+3d_j < \rho (p_{j+1})$. We may also require $\wt c_j > \rho(p_j)+3d_j$ and $\wt c_0> \wh c$. Define a new sequence
$$\begin{array}{lcl}
c_{2k}&=&\wt c_k,\\
c_{2k+1}&=&\rho(p_{k+1}) + d_{k+1},
\end{array}$$
where $k=0,1,2,\ldots\ $. If there are only finitely many critical points, choose the remainder of this sequence arbitrarily with $\lim_{j\to \infty} c_j =\infty$.

In the $j$-th stage of construction we assume inductively that we have a nonvanishing holomorphic 1-form $\omega_j$ in a neighborhood of $\{\rho \le c_j\}$, a holomorphic injective mapping $\alpha_j$ on a neighborhood of $\{\rho \le c_{j-1}\}$ and a number $\eps_j<\eps_{j-1}$ which satisfy the following:
\begin{enumerate}
\item[(1)] $\int_C \omega_j = \int_C\theta$ for each closed curve $C\subset \{\rho \le c_j\}$,
\item[(2)] $\omega_j = \alpha_j^*\omega_{j-1}$ and $||\alpha_j-id||<\eps_{j-1}$ on a domain of $\alpha_j$,
\item[(3)] for all $x,y\in \{\rho\le c_{j-1}\}$ the condition $d(x,y)<\eps_j$ implies
$$||\alpha_k\circ\alpha_{k+1}\circ \ldots \circ \alpha_j (x) - \alpha_k\circ\alpha_{k+1}\circ \ldots \circ \alpha_j (y)|| < \frac{\eps}{2^{j+1}}$$
for each number $k\in \{1,2,\ldots,j\}.$
\end{enumerate}
Since $\{\rho\le c_0\}$ is a noncritical strongly pseudoconvex extension of $\{\rho \le \wh c\}$, Proposition \ref{prop2} provides a closed holomorphic 1-form $\omega_0$ on a neighborhood of $\{\rho \le c_0\}$ such that $\omega_0 =\alpha_0^* \omega$ on a neighborhood of $\{\rho \le \wh c\}$, where $||\alpha_0-id||<\frac{\eps}{2}$. Choose $\eps_0$ such that for every $x,y\in \{\rho \le \wh c\}$ with $d(x,y)<\eps_0$ it follows that $||\alpha_0(x)-\alpha_0(y)||<\frac{\eps}{2}$. We now explain how to obtain $\omega_{j+1}$.

If there is no critical value on the interval $[c_j,c_{j+1}]$, we use Proposition \ref{prop2} to get $\omega_{j+1}$ and $\alpha_{j+1}$ which satisfy
\begin{equation}\label{prva}
\int_C \omega_{j+1} = \int_C\theta \textrm{ for each closed curve } C\subset \{\rho \le c_{j+1}\},
\end{equation}
and
\begin{equation}\label{druga}
\omega_{j+1} = \alpha_{j+1}^*\omega_{j},\ ||\alpha_{j+1}-id||<\eps_j \textrm{ on the domain of } \alpha_{j+1}.
\end{equation}
Now suppose that the interval $[c_j,c_{j+1}]$ contains a critical value. The construction in subsection \ref{critical} provides $\omega_{j+1}$ and $\alpha_{j+1}$ for which statements (\ref{prva}) and (\ref{druga}) hold. In both cases, we may (by uniform continuity) choose  $\eps_{j+1}<\eps_j$ in such a way that $d(x,y)<\eps_{j+1}$
implies
\begin{equation}\label{tretja}
||\alpha_k\circ\alpha_{k+1}\circ \ldots \circ \alpha_{j+1} (x) - \alpha_k\circ\alpha_{k+1}\circ \ldots \circ \alpha_{j+1} (y)|| < \frac{\eps}{2^{j+2}}
\end{equation}
for all $x,y\in \{\rho\le c_{j}\}$ and for each number $k\in \{1,2,\ldots,j+1\}.$

Define a global 1-form $\wt\omega$ as a limit $\wt \omega=\lim_{j\to \infty} \omega_j$. It remains to show that the sequence converges and $\omega$ satisfies the stated properties. Let $K$ be an arbitrary compact set in $X$. There is an index $j$ such that $K\subset \{\rho\le c_{j-1}\}$. Since $\omega|_K=\lim_{k\ge j} \omega_j|_K$, we have
$$\omega|_K=\lim_{k\to\infty} (\alpha_{j}\circ \ldots \circ \alpha_{j+k})^*\omega_j|_K.$$
Set $h_j:=\lim_{k\to\infty} (\alpha_{j}\circ \ldots \circ \alpha_{j+k})$. To see that $h_j$ is holomorphic, first choose any $\delta>0$. There is an index $k$ such that $2^{j+k-1} \cdot \delta > \eps$. For any positive integer $l$ we have
$$||\alpha_j\circ \ldots \circ \alpha_{j+k+l} - \alpha_j \circ \ldots \circ \alpha_{j+k}||_K=$$
$$\bigg|\bigg|\sum_{s=0}^{l-1}\left(\alpha_j\circ \ldots \circ \alpha_{j+k+s+1} - \alpha_j \circ \ldots \circ \alpha_{j+k+s}\right)\bigg|\bigg|_K\le$$
$$\le \sum_{s=0}^{l-1}||\alpha_j\circ \ldots \circ \alpha_{j+k+s}\circ \alpha_{j+k+s+1} - \alpha_j \circ \ldots \circ \alpha_{j+k+s}||_K.$$
Since by (\ref{druga}) we have $||\alpha_{j+k+s+1}-id||_K<\eps_{j+k+s}$, the condition (\ref{tretja}) implies that
$$||\alpha_j\circ \ldots \circ \alpha_{j+k+l} - \alpha_j \circ \ldots \circ \alpha_{j+k}||_K\le \sum_{s=0}^{l-1} \frac{\eps}{2^{j+k+s}}<\sum_{s=0}^{l-1} \frac{2^{j+k-1} \delta}{2^{j+k+s}} < \delta$$
for any positive integer $l$. Since each $\alpha_k$ is injective, it follows by chosing the norms precise enough that $h_j$ is holomorphic and injective. Therefore $\omega$ is a closed holomorphic 1-form on $X$ that is without zeros. For any closed curve $C\subset K$ we have
$$\int_C \omega =\int _C h^*\omega_j=\int_{h(C)}\omega_j = \int _C \theta$$
by Stokes theorem. Since $K$ was chosen arbitrarily, the conclusion holds for any closed curve $C\subset X$, which implies that $\omega$ is in the same cohomology class as $\theta$.

\section{Linearly independent 1-forms}


In this section we prove the following generalization of Theorem \ref{glavni} and Theorem II in \cite{For1}.

\begin{theorem}
Let $X^n$ be a Stein manifold whose holomorphic cotangent bundle $T^*X$ admits a trivial complex subbundle of rank $q$ for some $1\le q<n$. Given closed 1-forms $\theta_1, \ldots, \theta_q$ on $X$ there exist closed holomorphic 1-forms $\omega_1, \ldots, \omega_q$ satisfying
$$[\omega_j]=[\theta_j] \in H^1(X,\C) \quad \textrm{for each $j=1,\ldots, q$}$$
and
$$\omega_1 \wedge \omega_2 \wedge \ldots \wedge \omega_q |_x \ne 0\quad \textrm{for all $x\in X$}.$$
\end{theorem}

\begin{remark}
$T^*X$ always contains a trivial subbundle of rank $q= [\frac{n+1}{2}]$ and this value is optimal in general (see \cite[p.~714]{For-1} and \cite[Proposition 3]{For0}).
\end{remark}

\begin{proof}{}
We choose continuous $(1,0)$-forms $\omega_1',\ldots, \omega_q'$ such that
$$\omega_1'\wedge \ldots \wedge \omega_q'|_x\ne 0$$
for all $x\in X$. These forms span a trivial subbundle of $T^*X$. We construct pointwise linearly independent holomorphic 1-forms $\omega_1, \ldots, \omega_q$ such that $[\omega_j]=[\theta_j] \in H^1(X,\C)$ for each $j=1,\ldots, q$, where the collection $(\omega_1, \ldots, \omega_q)$ is homotopic to $(\omega'_1, \ldots, \omega_q')$ through the homotopy of $q$-tuples of independent sections of $T^*X$.  The construction is similar to the proof of Theorem \ref{applic} and we briefly illustrate the necessary changes in the noncritical and the critical case.\medskip

\paragraph{(I.) The noncritical case.} Choose a local representation of a bump as in subsection \ref{noncritical} and represent each $\omega_j$ by a noncritical holomorphic function $f_j$ on a neighborhood of $\wt K$. Approximate a submersion $f=(f_1, \ldots, f_q)$ by a holomorphic submersion $g\colon V'\to \C^q$, $V'\supset \wt L$, using Proposition 3.3.~in \cite{For1}. By Lemma 5.1.~in \cite{For1} there is an injective holomorphic map $\gamma \colon W\to X$, $W\subset V'$, satisfying $f=g\circ \gamma$ on $W$ and $||\gamma - id||_{W} < \eps$. Using Theorem \ref{split} we obtain $\gamma = \beta \circ \alpha^{-1}$ with $\alpha$ and $\beta$ close to identity. Finally define $\wt\omega_j$ as
$$\wt \omega_j =\left\{\begin{array}{c}
\alpha ^* \omega_j,\\
d(g_j\circ \beta).
\end{array}\right.$$
By construction $\wt \omega_j$ has the same periods as $\theta_j$ and no zeros. Since
$$\wt\omega_1\wedge\ldots\wedge\wt\omega_q=\left\{
\begin{array}{cc}
\alpha^*(\omega_1\wedge \ldots \wedge \omega_q)& \textrm{on } A',\\
\beta^*(dg_1\wedge \ldots \wedge dg_q) & \textrm{on } B',
\end{array}\right.$$
the 1-forms $\wt\omega_1, \ldots, \wt \omega_q$ are pointwise linearly independent.

\bigskip

\paragraph{(II.) The critical case.} The pointwise independency of 1-forms $\omega_1$, $\ldots$, $\omega_q$ is invariant under small perturbations, hence there is an $\eps_\alpha>0$ such that the conditions $||\alpha_j-id||_{A'}<\eps_\alpha$ for $j=1,\ldots,q$ imply that
$$\alpha_1^* \omega_1 \wedge \ldots \wedge\alpha_q^* \omega_q|_x\ne 0 \quad \textrm{for all } x\in A'.$$

\paragraph{(II.1) The Morse index is at least 2.} We use the notations as in subsection \ref{critical}.  There is a neighborhood $V$ of $bE$ in $\{\rho\le -\frac{c}{2}\}$ and a holomorphic submersion $f\colon V\to \C^q$ such that on $V$ we have local representations $\omega_j = df_j$ for $j=1, \ldots, q$. 

By Lemma 6.4.~in \cite{For1} there is a constant $c'>c$ close to $c$ such that $f$ and it's partial derivatives restricted to $\{\rho \le -c'\}$ extend smoothly to $\wt f$ on $(\{\rho \le -c'\} \cap V)\cup E$ where the Jacobian matrix $J(\wt f)=(\partial\wt f_j/\partial z_l)$ has complex rank $q$ at each point of $E$.

By Lemma 6.6.~in \cite{For1} we get for every $\delta>0$ a holomorphic submersion $g\colon V'\to \C^q$ such that $|\wt f-g|_K<\delta$, $|d\wt f-dg|_E<\delta$. A small perturbation of each element of the set $\{dg_1,\ldots , dg_q\}$ leaves pointwise independency invariant, thus there is $\eps_\beta>0$ such that  $d(g_1\circ \beta_1),\ldots, d(g_q\circ \beta_q)$ are pointwise independent if $||\beta_j -id||< \eps_\beta$.

We use Lemma 5.1.~in \cite{For1} and Theorem \ref{split} to get $f_j=g_j\circ\gamma_j$ and split $\gamma_j=\beta_j\circ \alpha_j^{-1}$. This can be done in such a way that $||\alpha_j-id||<\eps_\alpha$ and $||\beta_j-id||<\eps_\beta$ for each $j$. The obtained 1-forms
$$\wt\omega_j=\left\{\begin{array}{cc}
\alpha_j^*\omega_j & \textrm{on } A',\\
d(g_j\circ \beta_j)& \textrm{on } B',
\end{array}\right.$$
are then pointwise independent with prescribed periods and no zeros.\medskip

\paragraph{(II.2) Morse index of a critical point is 1.} In this case the neighborhood $V$ of $bE$, where $df_j=\omega_j$, can be chosen as a disjoint union of two open sets $V_a$ and $V_b$. We may further assume that
$$f_j(b)-f_j(a)=\int_{C_0} \theta_j - \int_{C_0\bs E}\omega_j.$$

Let $f_j^{t_j}\colon V\to \C$ be defined as $f_j^{t_j}|_{V_a}=f_j$ and $f_j^{t_j}|_{V_b}=f_j^{t_j}+t_j$. Denote $f_{t_1\ldots t_q} = (f_1^{t_1}, \ldots , f_q^{t_q})\colon V\to\C^q$. Lemma 6.4.~and Lemma 6.6.~in \cite{For1} can be used for a family of submersions $\{f_{t_1\ldots t_q} \}$ holomorphically depending on parameters $(t_1,\ldots, t_q)\in \Delta\times\ldots\times\Delta$ to get families $\wt f_{t_1\ldots t_q}$ and $g_{t_1\ldots t_q}$. Again we may find such number $\eps_\beta>0$  that the condition $||\beta_j^{t_j} - id||<\eps_\beta$  implies the pointwise independency of $d(g_1^{t_1}\circ \beta_1^{t_1}), \ldots, d(g_q^{t_q}\circ \beta_q^{t_q})$. 

By Lemma 5.1.~in \cite{For1} we get for each $j=1,\ldots, q$ maps $\gamma_j ^{t_j}$ such that $f_j^{t_j} = g_j^{t_j} \circ \gamma_j^{t_j}$. Splitting $\gamma_j(x,t)=(\gamma_j^{t}(x), t)$ with respect to $A\times \Delta$ and $B\times \Delta$ we get $\gamma_j=\beta_j\circ \alpha_j^{-1}$ where
$$\alpha_j(x,t)=(\alpha_j^t(x),t),\ \beta_j(x,t)=(\beta_j^t(x),t), \
||\alpha_j -id||<\eps_\alpha,\ ||\beta_j-id||<\eps_\beta.$$
We finally define
$$\wt\omega_j^{t_j} = \left\{\begin{array}{c}
(\alpha_j^{t_j})^*\omega_j,\\
d(g_j^{t_j}\circ \beta_j^{t_j}).
\end{array}\right.$$
Using Rouch\'e's theorem as at the end of subsection \ref{critical}, we find $\wt t_j\in\Delta$ in such a way that
$$\int_{C_0}\wt \omega_j^{\wt t_j} = \int _{C_0} \theta_j.$$
By construction the 1-forms $\wt\omega_1^{\wt t_1},\ldots, \wt\omega_q^{\wt t_q}$ have no zeros and are pointwise independent.
\qed\end{proof}

\section{Algebraic example}

Due to Forstnerič \cite{For1} each Stein manifold $X$ admits a noncritical holomorphic function (for an open Riemann surface see also \cite{GunNar}). When $X\subset \C^N$ is an affine algebraic manifold, a natural question is whether there exists a holomorphic polynomial $P$ on $\C^N$ whose restriction to $X$ is noncritical to $X$. The following counter example was provided by R.~Narasimhan (private communication).

Let $t\in \C$ and let $\Gamma\subset \C$ be a lattice of rank two. Let $x(t)$,  $y(t)$ be $\Gamma$-periodic meromorphic functions (elliptic functions) satisfying an equation
$$y^2 = x^3+Ax+B.$$
Denote by $C$ the corresponding curve in $\C^2$.

\begin{proposition}
For every polynomial $P\in \C[x,y]$ the restriction
$$P|_C\colon C\to\C$$
has at least one critical point on $C$.
\end{proposition}

\begin{proof}
Let $p(t)=P(x(t),y(t))$ for $t\in \C\bs \Gamma$. Clearly $p$ and its derivative $\dot p$ are $\Gamma$-periodic functions, and all points of the lattice are singularities of the same kind due to periodicity. Since $P$ is a polynomial and $x(t)$, $y(t)$ are meromorphic, they cannot be essential singularities.

If $\dot p$ is regular at points of $\Gamma$, it follows that it is bounded on $\C$ and hence constant, thus $\dot p=c$. Hence $p(t)=ct+b$ for some $b\in\C$ but this function fails to be periodic. Thus $\dot p$ must have a pole at every point of $\Gamma$ and hence $1/\dot p$ has zeros at these points.

If $\dot p$ has no zeros on $\C$, it follows that $1/\dot p$ is a bounded function on $\C$ which is a contradiction. This means that $\dot p(t_0)=0$ for some $t_0\in\C\bs\Gamma$, which means that $P|_C$ has a critical point at $(x(t_0), y(t_0))\in C$.
\qed\end{proof}

\small{\noindent {\bf Acknowledgements.} 
I am very grateful to F.~Forstneri\v c for helpful discussions. I wish to thank R. Narasimhan for suggesting the example in section 4 and for proposing that Theorem \ref{glavni} may hold on open Riemann surfaces (private communication, 2004).

Research on this work was supported by a fellowship from the Agency for Research and Development (ARRS)
of the Republic of Slovenia (contract no.~1000-05-310002).


\end{document}